\documentclass[a4paper,12pt]{article}
\newcommand{\stoppage}{\newpage{\pagestyle{empty}\cleardoublepage}}

\newtheorem{theorem}{Theorem}[section]
\newtheorem{lemma}[theorem]{Lemma}
\newtheorem{corollary}[theorem]{Corollary}
\newenvironment{proof}{\begin{trivlist}\item[]\textbf{Proof}}%
                       {\setqed\end{trivlist}}
\newenvironment{remark}{\begin{trivlist}\item[]\textbf{Remark}}{\end{trivlist}}
\newenvironment{acknow}{\begin{trivlist}\item[]\textbf{Acknowledgements}}
 {\end{trivlist}}

\usepackage{epsfig}
\usepackage{latexsym}
\usepackage{amssymb}

\newcommand{\epsdir}{eps}
\newcommand{\setqed}{\hfill\ensuremath{\diamond}}

\newcommand{\sa}{{\cal C}}
\newcommand{\sap}{{\cal C'}}
\DeclareSymbolFont{AMSb}{U}{msb}{m}{n}
\DeclareMathSymbol{\myZ}{\mathalpha}{AMSb}{'132}

\newcommand{\id}{\mathrm{id}}

\newcommand{\sothat}{\;|\;}
\newcommand{\tensor}{\otimes}
\newcommand{\done}{\Box}
\newcommand{\bigdone}{\Box}

\newcommand{\qi}[1]{[#1]}

\begin{document}

\begin{center}
{\LARGE
Idempotents of the Hecke algebra become \\  Schur functions in
 the skein of the annulus}
\vspace{5mm}
\\ {\large Sascha G. Lukac}
\vspace{2mm}
\\ {\ttfamily\large lukac@liv.ac.uk}
\vspace{3mm}
\\ {\large Department of Mathematical Sciences, University of Liverpool, 
    Peach St, Liverpool L69 7ZL, England}
\end{center}
%%%%%

\begin{abstract}
The Hecke algebra $H_n$ contains well known idempotents $E_{\lambda}$
which are indexed by Young diagrams with $n$ cells. They were originally
described by Gyoja \cite{Gyoja}. 

A skein theoretical description of $E_{\lambda}$
was given by Aiston and Morton \cite{AisMor}. The closure of $E_{\lambda}$ becomes
an element $Q_{\lambda}$ of the skein of the annulus. In this skein, they are known to 
obey the same multiplication rule as the symmetric Schur functions $s_{\lambda}$
as stated in theorem \ref{T-ringhom}.
But previous proofs of this fact as in \cite{Aiston} used results about quantum groups
which were far beyond the scope of skein theory.

Our elementary proof was motivated by \cite{Kawagoe} and uses only skein theory and
basic algebra.
\end{abstract}

\section{The skein of a planar surface}\label{sec-skeinpl}

\begin{figure}[!h]
  \begin{center}
    \setlength{\unitlength}{1mm}
    \begin{picture}(0,0)
    \put(11,7){$-$}
    \put(28,7){$=\ \ (s-s^{-1})\cdot$}
    \put(66,7){,}
    \put(96,7){$=\ v^{-1}\cdot$}
    \end{picture}
     \epsfig{file=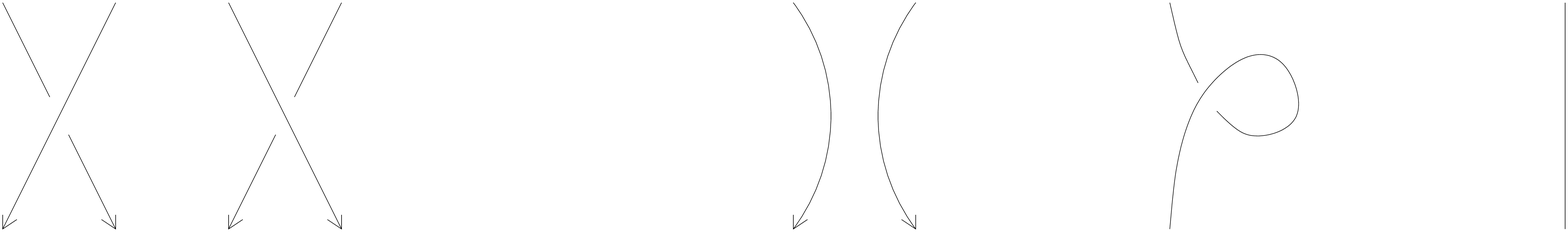,width=110\unitlength}
     \caption{\label{fig-skeinrel}Skein relations.}
  \end{center}
\end{figure}

We consider a planar surface $F$ with designated $n$ incoming and $n$ outgoing
boundary points for some integer $n\geq 0$.
The Homfly skein ${\cal S}(F)$ is defined as the module of linear combinations
of oriented 
tangles in $F$ quotiented by regular isotopy (i.e. Reidemeister moves II and III),
the two local skein relations in figure \ref{fig-skeinrel}
where $v$ and $s$ are variables,
and the relation that a disjoint simple closed curve can be removed from a diagram at the
expense of multiplication with the scalar $\delta=\frac{v^{-1}-v}{s-s^{-1}}$.
The last relation is a consequence of the other relations unless
the remaining diagram is the empty diagram.

The scalars can be chosen as $\myZ[v^{\pm 1},s^{\pm 1}]$ 
with powers of $s-s^{-1}$ in the denominators. If necessary, we can choose
the scalars to be the field of rational functions in $v$ and $s$.
This will occur in section \ref{sec-hecke} when we divide a quasi-idempotent by
a scalar to turn it into an idempotent.

For the skein of any planar surface $F$ there is
a map $\rho:S(F)\to S(F)$ which is induced by switching all crossings
and replacing $s$ by $s^{-1}$ and $v$ by $v^{-1}$.
The skein relations are preserved by this map. We remark that
$\rho$ is not linear.

\section{Young diagrams and symmetric functions}
\subsection{Young diagrams}

We consider partitions of non-negative integers where the
summands are in non-strictly decreasing order. We consider partitions
to be equal if they differ only by a number of zero-summands at the end.

A \emph{Young diagram} is a graphical depiction of a partition.
A partition $\lambda=(\lambda_1,\lambda_2,\ldots,\lambda_k)$
becomes a flush-left arrangement of square cells with $\lambda_1$ cells
in the first row, $\lambda_2$ cells in the second row, ... and $\lambda_k$
cells in the $k$-th row. 

We define $\lambda\subset\mu$ for partitions
$\lambda=(\lambda_1,\ldots,\lambda_k)$ and $\mu=(\mu_1,\ldots,\mu_m)$
if $\lambda_i\leq \mu_i$ for all $1\leq i\leq k$ (where $\mu_i=0$ for $m<i\leq k$).
This relation coincides in the graphical depiction with $\lambda$
being a subset of $\mu$.

We define $|\lambda|$ to be the number of cells of $\lambda$, i.e.
for $\lambda=(\lambda_1,\ldots,\lambda_k)$ we have $|\lambda|=
\lambda_1+\cdots+\lambda_k$.

\subsection{Symmetric functions}

Following Macdonald \cite{Macdonald}, we consider the polynomial ring in variables
$x_1,\ldots,x_n$. For a partition $\lambda=(\lambda_1,\ldots,\lambda_k)$
with $k\leq n$ we define
\[
  s_{\lambda}(x_1,\ldots,x_n)=\frac{\det(x_i^{\lambda_j+n-j})_{1\leq i,j\leq n}}
          {\det(x_i^{n-j})_{1\leq i,j\leq n}}.
\]
This quotient is a symmetric polynomial in $x_1,\ldots,x_n$.
In the inverse limit of the polynomial rings $\myZ[x_1,\ldots,x_n]$
for $n\to\infty$ we can define an element $s_{\lambda}$,
called the \emph{$\lambda$-Schur function} for any Young diagram $\lambda$.
This element specializes to $s_{\lambda}(x_1,\ldots,x_n)$ when we set 
$x_i=0$ for all $i\geq n+1$ provided that $\lambda$ has at most $n$ non-zero summands.

The most basic Young diagrams are on the one hand a row diagram with, say, $n$ cells
and on the other hand a column diagram with, say, $k$ cells. One can check that
$s_{(n)}$ is the $n$-th complete symmetric function $h_n$ and that $s_{(1^k)}$
is the $k$-th elementary symmetric function $e_k$.

\section{The skein $R_n^n$ of a disk with $2n$ boundary points}\label{sec-hecke}

We consider a rectangle $F$ with $n$ outgoing points at the top and $n$ incoming
points at the bottom, $n\geq 0$.
We denote the Homfly skein of this surface by $R_n^n$. It is known to be isomorphic
to the Hecke algebra $H_n$.

$R_n^n$ gets an algebra structure by defining $D\cdot E$ as putting $D$ above $E$
so that the incoming points of $D$ match with the outgoing points of $E$.
Clearly, $R_n^n$ is generated as an algebra by the elementary braids $\sigma_1,\ldots,
\sigma_{n-1}$ where $\sigma_i$ agrees with the identity braid up to 
a single positive crossing between the arcs which connect the boundary points $i$ and 
$i+1$ (numbered from the left).

As a module, $R_n^n$ is linearly spanned by the positive permutation braids
$\omega_{\pi}$ where $\pi$ is a permutation on $n$ letters (see \cite{Nato} 
for details).

$R_n^n$ has a well known quasi-idempotent $a_n$ which is given
in terms of the positive permutation braids $\omega_{\pi}$ 
on $n$ strings as
\[
  a_n=\sum_{\pi\in S_n}s^{l(\pi)}\omega_{\pi}
\]
where $l(\pi)$ is the writhe of the braid $\omega_{\pi}$.
This element `swallows' an elementary braid $\sigma_i$ at the
expense of a scalar $s$ (see \cite{Nato}).
\begin{lemma}\label{L-ans}
$a_n\sigma_i=\sigma_ia_n=sa_n$ for any $1\leq i\leq n-1$.
\end{lemma}
This behaviour implies that $a_n$ is central in $R_n^n$ because
the elementary braids $\sigma_1,\ldots,\sigma_{n-1}$ generate $R_n^n$. 
Lemma \ref{L-ans} enables an inductive proof that
$a_na_n=\alpha_na_n$ where the scalar $\alpha_n$ is given by
\[
  \alpha_n=s^{\frac{n(n-1)}{2}}\qi{n}\qi{n-1}
   \cdots\qi{1}
\]
where $\qi{i}=\frac{s^i-s^{-i}}{s-s^{-1}}$ for any integer $i\geq 0$.

We recall the map $\rho$ from section \ref{sec-skeinpl}.
\begin{lemma}\label{L-alpn}
  $\rho\left(\frac{1}{\alpha_n}a_n\right)=\frac{1}{\alpha_n}a_n$
for any $n\geq 1$.
\end{lemma}
\begin{proof}
We have
$\rho(ab)=\rho(a)\rho(b)$ for any elements $a,b\in R_n^n$.
This multiplicativity together with lemma \ref{L-ans} implies that
\[ 
  \sigma_i\rho(a_n)=\rho(\sigma_i^{-1})\rho(a_n)
  = \rho(\sigma_i^{-1}a_n)=\rho(s^{-1}a_n)=s\rho(a_n)
\]
in $R_n^n$. Hence, $a_n\rho(a_n)=\alpha_n\rho(a_n)$. When we apply
the map $\rho$ to this equation we get $\rho(a_n)a_n=\rho(\alpha_n)a_n$.
The left hand sides of the previous two equations are equal because
$a_n$ is central. Hence, $\alpha_n\rho(a_n)=\rho(\alpha_n)a_n$ and thus
$\rho\left(\frac{1}{\alpha_n}a_n\right)=\frac{1}{\alpha_n}a_n$.
\end{proof}

Aiston and Morton \cite{AisMor} constructed idempotents $E_{\lambda}$
in $R_n^n$ for any Young diagram $\lambda$ with $n$ cells. To do this,
they considered a three-dimensional version of $R_n^n$ where the distinguished
boundary points are lined up along the cells of the Young diagram $\lambda$.
In this setting they arranged several copies of idempotents $a_i$
and varaints $b_i$ of them along the rows and columns of $\lambda$ to get
a quasi-idempotent. This element which is denoted $e_{\lambda}$ in the
equivalent $R_n^n$-setting satisfies $e_{\lambda}e_{\lambda}=
\alpha_{\lambda}e_{\lambda}$. It is sufficient to prove that $\alpha_{\lambda}\neq 0$
for the specialization $v=s=1$ which turns $R_n^n$ into the
symmetric-group algebra as explained in \cite{Aiston}.

We consider the rational functions in $v$ and $s$ as the set of scalars for the
skein $R_n^n$. Then,
the element $E_{\lambda}=\frac{1}{\alpha_{\lambda}}e_{\lambda}$
is an idempotent of $R_n^n$.

\section{The skein $\sa$ of the annulus}\label{sec-skeinan}
\begin{figure}
  \begin{center}
    \setlength{\unitlength}{1mm}
    \begin{picture}(0,0)
    \put(36,15){$\cdot$}
    \put(75,15){=}
    \end{picture}
     \epsfig{file=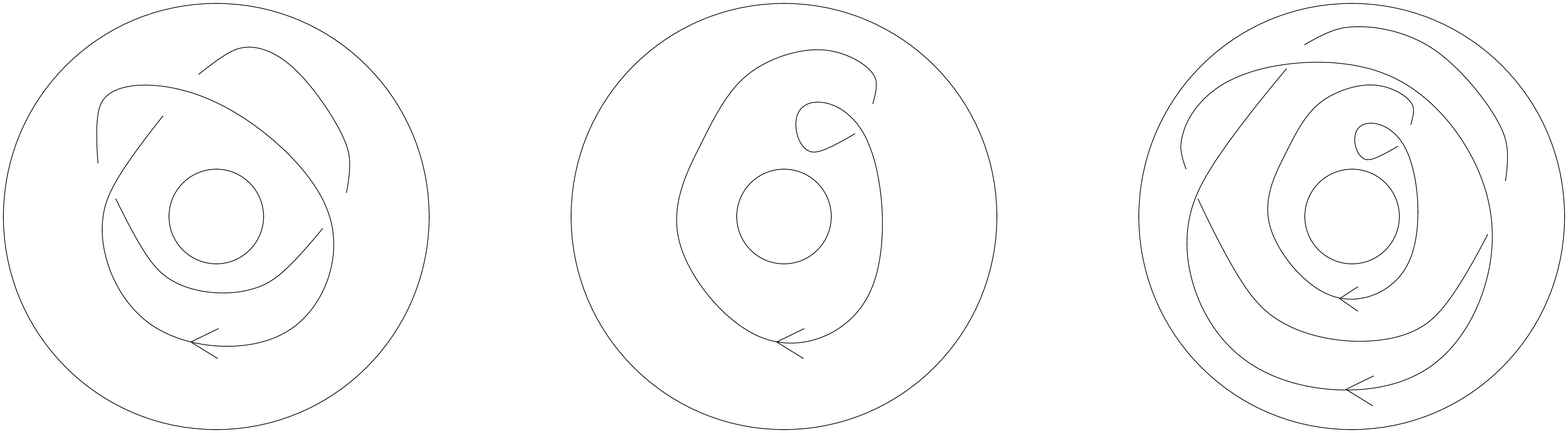,width=110\unitlength}
     \caption{\label{F-anmult} The multiplication in the skein of
               the annulus $\sa$.}
  \end{center}
\end{figure}

The skein of the annulus shall be denoted by $\sa$.
Let $D_1$ and $D_2$
be two diagrams in the annulus $S^1\times [0,1]$.
We can bring $D_1$ into $S^1\times [0,1/2)$, and $D_2$
into $S^1\times (1/2,1]$ by regular isotopies. Then the
product of $D_1$ and $D_2$ is defined as the diagram
$D_1\cup D_2$.
An example is shown in figure \ref{F-anmult}.
 The product is commutative since $D_1D_2$
and $D_2D_1$ differ by regular isotopy.
The empty diagram is the identity.

\begin{figure}
  \begin{center}
    \setlength{\unitlength}{1mm}
    \begin{picture}(0,0)
    \end{picture}
     \epsfig{file=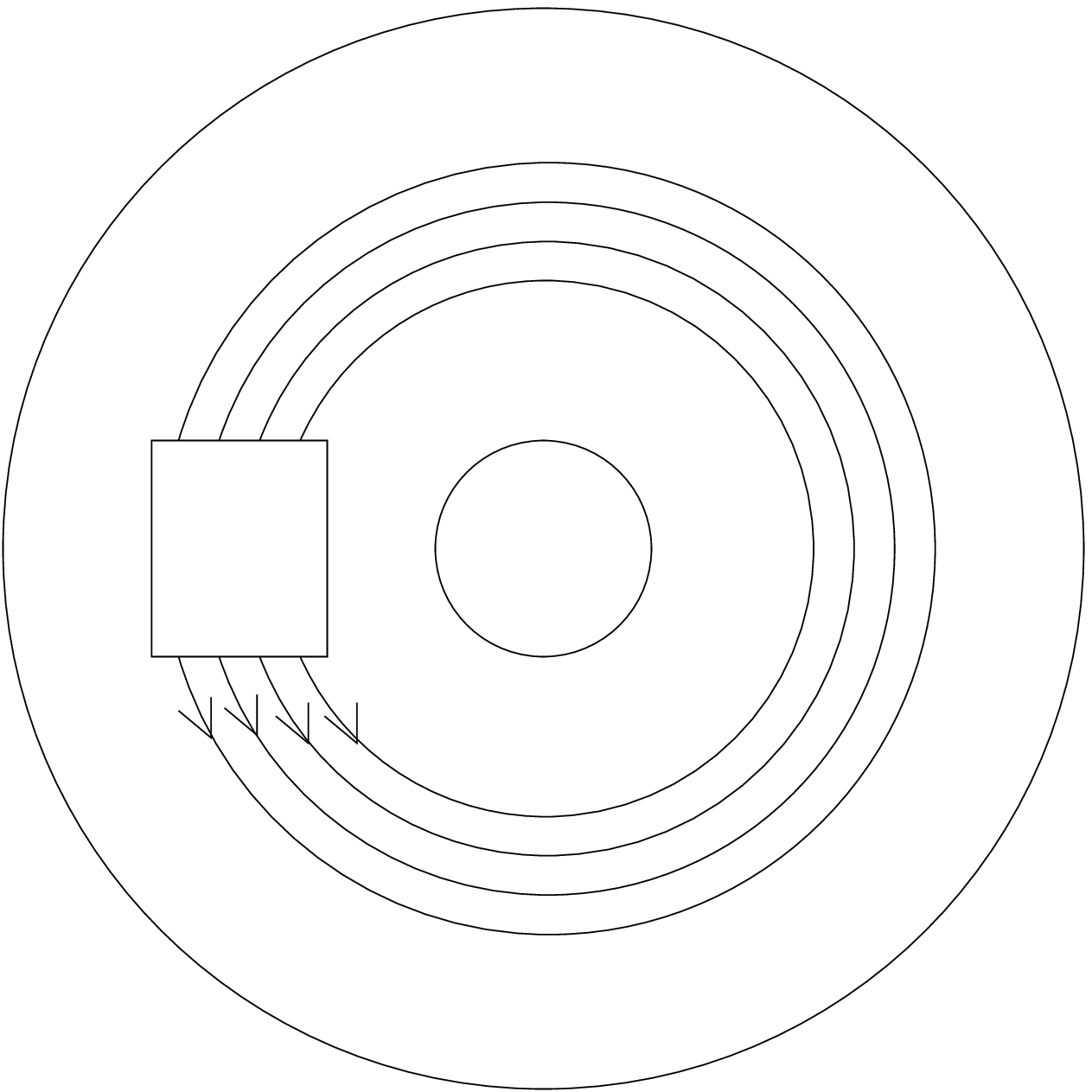,width=40\unitlength}
     \caption{\label{F-cls} The additional arcs for the closure map $\Delta_n:R_n^n\to\sa$.}
  \end{center}
\end{figure}

Figure \ref{F-cls} depicts an annulus
with a set of $n$ oriented arcs.
A rectangle is removed from the annulus so that we can insert
a diagram from $R_n^n$ with matching orientations of the arcs.
This factors to a map $\Delta_n:R_n^n\to\sa$.
This is a special case of a \emph{wiring}. 
We abbreviate $\Delta_n$ by $\Delta$ if the index $n$ is obvious from
the context.

We denote by $\sa_n$ the image of $R_n^n$ under the closure map
which is a submodule of $\sa$.
 By $\sa_+$
 we denote the submodule of $\sa$ spanned
by all $\sa_0,\sa_1,\ldots$.

We define $Q_{\lambda}$
to be the closure of the idempotent $E_{\lambda}$ of $R_n^n$
where $n$ is the number of cells of $\lambda$,
\[
  Q_{\lambda}=\Delta_n(E_{\lambda})\in \sa.
\]
By looking at the limit $s\to 1$ it is straightforward to verify
that the Homfly polynomial of $Q_{\lambda}$ as an element of the plane
is non-zero and therefore $Q_{\lambda}$ is non-zero in $\sa$ (see \cite{Lukac}).

\begin{figure}
  \begin{center}
    \setlength{\unitlength}{1mm}
    \begin{picture}(0,0)
    \put(44,17){$\stackrel{\Gamma}{\longrightarrow}$}
    \end{picture}
     \epsfig{file=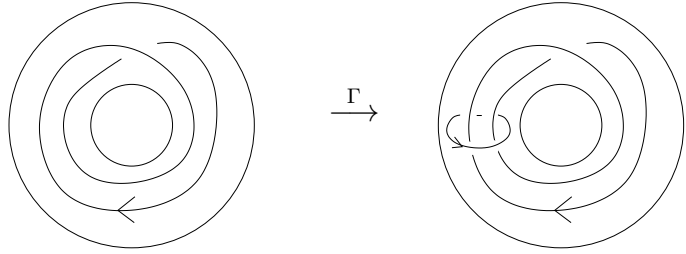,width=90\unitlength}
     \caption{\label{F-encirc} Encircling a diagram in the annulus.}
  \end{center}
\end{figure}

We define a linear map $\Gamma$
 from $\sa_+$ to $\sa_+$ that is
the encircling of a diagram in the annulus by a single loop with
a specified orientation
as shown in figure \ref{F-encirc}.

In $R_n^n$, one can write the identity braid on $n$ strings with
an encircling loop as the linear combination of so-called Jucys-Murphy
operators as explained in \cite{Morton}.
The skein calculations in \cite{AisMor} then 
immediately imply the following lemma.
\begin{lemma}\label{L-qeigen}
For any Young diagram $\lambda$ we have
$\Gamma(Q_{\lambda})=c_{\lambda}Q_{\lambda}$
in $\sa_+$ with the scalar
\[
  c_{\lambda}=\frac{v^{-1}-v}{s-s^{-1}}+vs^{-1}
   \sum_{k=1}^{l(\lambda)}(s^{2(k-\lambda_k)}-s^{2k}).
\]
\end{lemma}
\begin{remark}
It is easy to confirm that
the eigenvalues $c_{\lambda}$ are pairwise different
for all Young diagrams $\lambda$.
\end{remark}

\begin{lemma}\label{L-dimcn}
The set $\{Q_{\lambda}\sothat \lambda \mbox{ has } n \mbox{ cells }\}$
 is a basis for
$\sa_n$ for any $n\geq 0$.
\end{lemma}
\begin{proof}
We define $X_i^+\in \sa_i$ as
the closure of the braid $\sigma_{i-1}\sigma_{i-2}\cdots\sigma_1\in R_i^i$.
The \emph{weighted degree} of a monomial
$(X_{i_1}^+)^{j_1}\cdots(X_{i_k}^+)^{j_k}$
is defined as $i_1j_1+\cdots + i_k j_k$.

Any diagram in the annulus can be written inductively via the skein relations
as a linear combination of totally descending curves. This means that
$\sa_n$ is spanned linearly by the monomials in the $X_i^+$ of weighted
degree $n$. Hence, the dimension of $\sa_n$ is at most $\pi(n)$, by which
we denote the number of partitions of $n$.

On the other hand, $Q_{\lambda}$ lies in $\sa_n$ provided that
the Young diagram $\lambda$ has $n$ cells. These elements are
non-zero and they are linearly independent since they have pairwise
different eigenvalues under the map $\Gamma$ by lemma
\ref{L-qeigen}.
By definition, there are $\pi(n)$ Young diagrams with $n$ cells.
Since the dimension of $\sa_n$ is at most $\pi(n)$ by the above argument,
the dimension of $\sa_n$ is exactly $\pi(n)$ and the set of $Q_{\lambda}$
where $\lambda$ has $n$ cells is a basis.
\end{proof}
\begin{remark}
Lemma \ref{L-dimcn} implies that the elements $X_i^+, i\geq 0$,
generate $\sa_+$ freely as a commutative algebra.
This is a special case of Turaev's result \cite{TuraevC}
that $\sa$ is generated freely as a commutative algebra by
$X_i^+,X_i^-$ for $i\geq 0$ where $X_i^-$ derives from
$X_i^+$ by reversing the orientation.
\end{remark}

\section{The variant skein $\sap$ of the annulus}

The skein of the annulus $\sap$ with two
boundary points has been considered e.g. in \cite{Kawagoe} and \cite{Gilmer}.
We are using the version introduced by Morton in \cite{Morton}
where the two boundary points lie on different boundary components
of the annulus. This version has the benefit of 
a commutative multiplication.

\begin{figure}
  \setlength{\unitlength}{1mm}
  \begin{center}
   \begin{picture}(0,0)
     \put(40,16){$\cdot$}
     \put(86,16){$=$}
   \end{picture}
     \epsfig{file=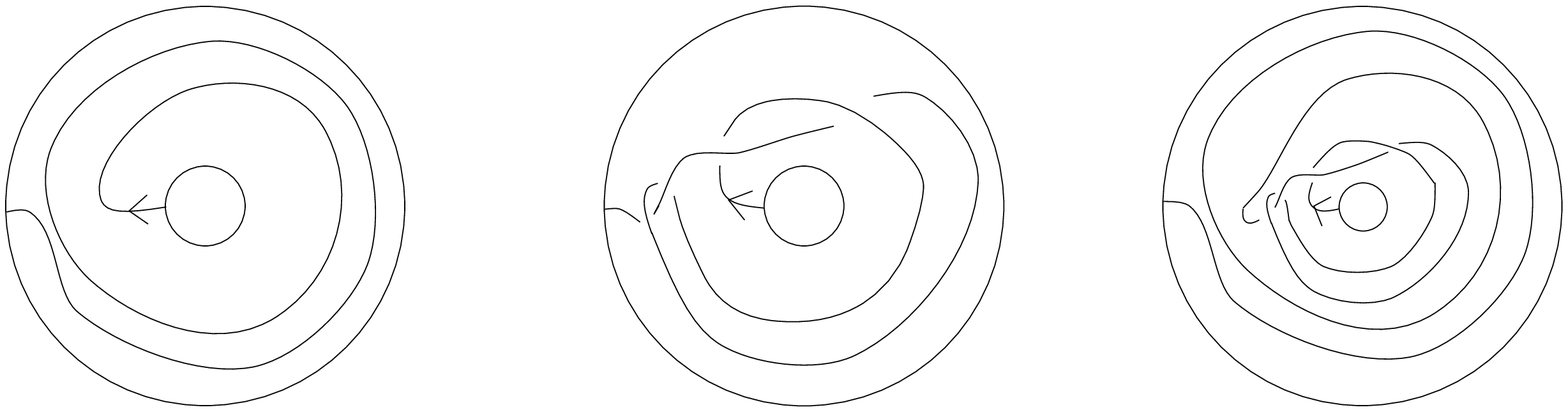,width=125\unitlength}
     \caption{\label{F-multcprime}The multiplication in $\sap$.}
  \end{center}
\end{figure}

\begin{figure}
 \setlength{\unitlength}{1mm}
 \begin{minipage}[t]{60\unitlength}
  \begin{center}
    \setlength{\unitlength}{1mm}
     \epsfig{file=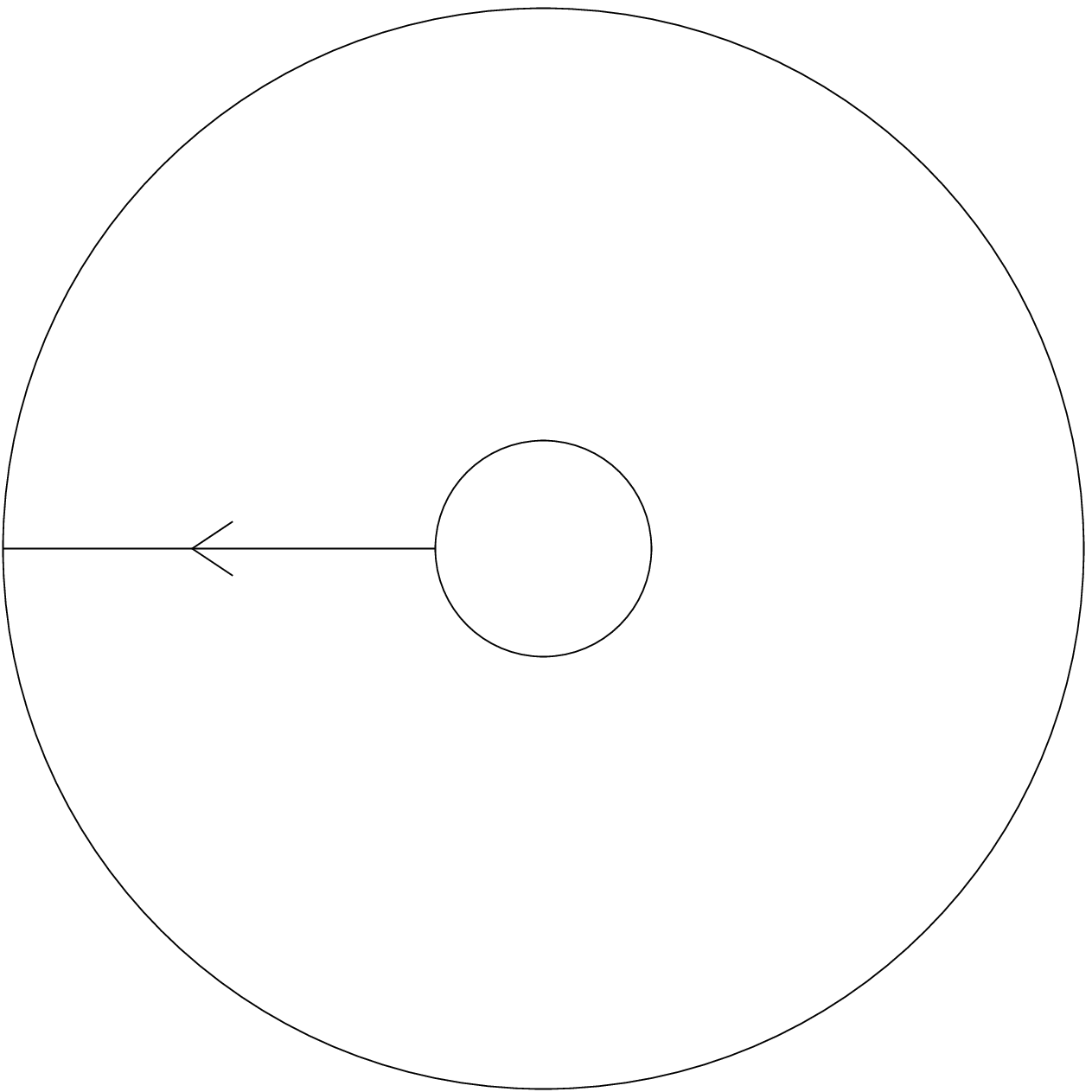,width=30\unitlength}
     \caption{\label{F-eiden}The identity $e$ in $\sap$.}
  \end{center}
 \end{minipage}
  \hfill
 \begin{minipage}[t]{70\unitlength}
  \begin{center}
    \setlength{\unitlength}{1mm}
     \epsfig{file=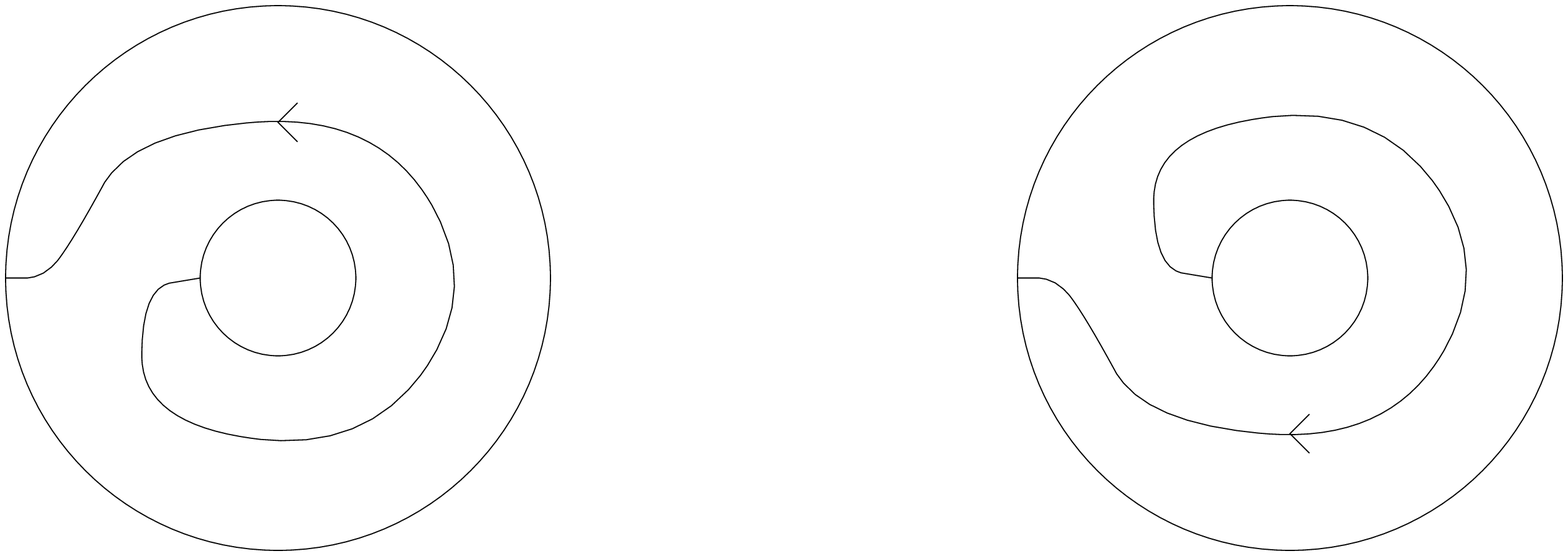,width=50\unitlength}
     \caption{\label{F-arca}The arc $a$ (at the left) 
                 and its inverse $a^{-1}$ (at the right).}
  \end{center}
 \end{minipage}
\end{figure}

We equip the annulus with two distinguished boundary points. An output
point on the inner circle, and an input point on the outer circle.
The resulting skein $\sap$ becomes an algebra by defining the product $A\cdot B$
of diagrams $A$ and $B$ as putting the outward circle of $B$ next to the inward
circle of $A$ so that the two involved distinguished boundary points become
a single point in the interior of the new annulus.
An example for the multiplication is depicted in figure \ref{F-multcprime}.
The single straight arc $e$ connecting the two marked points is the identity element,
as shown in figure \ref{F-eiden}.

Any arc which connects the two marked points without crossings is regularly
isotopic to a power of the diagram $a$ which is depicted in figure \ref{F-arca}.

\begin{figure}
  \begin{center}
    \setlength{\unitlength}{1mm}
    \begin{picture}(0,0)
    \put(5,17){$D$}
    \put(40,17){$D$}
    \end{picture}
   \epsfig{file=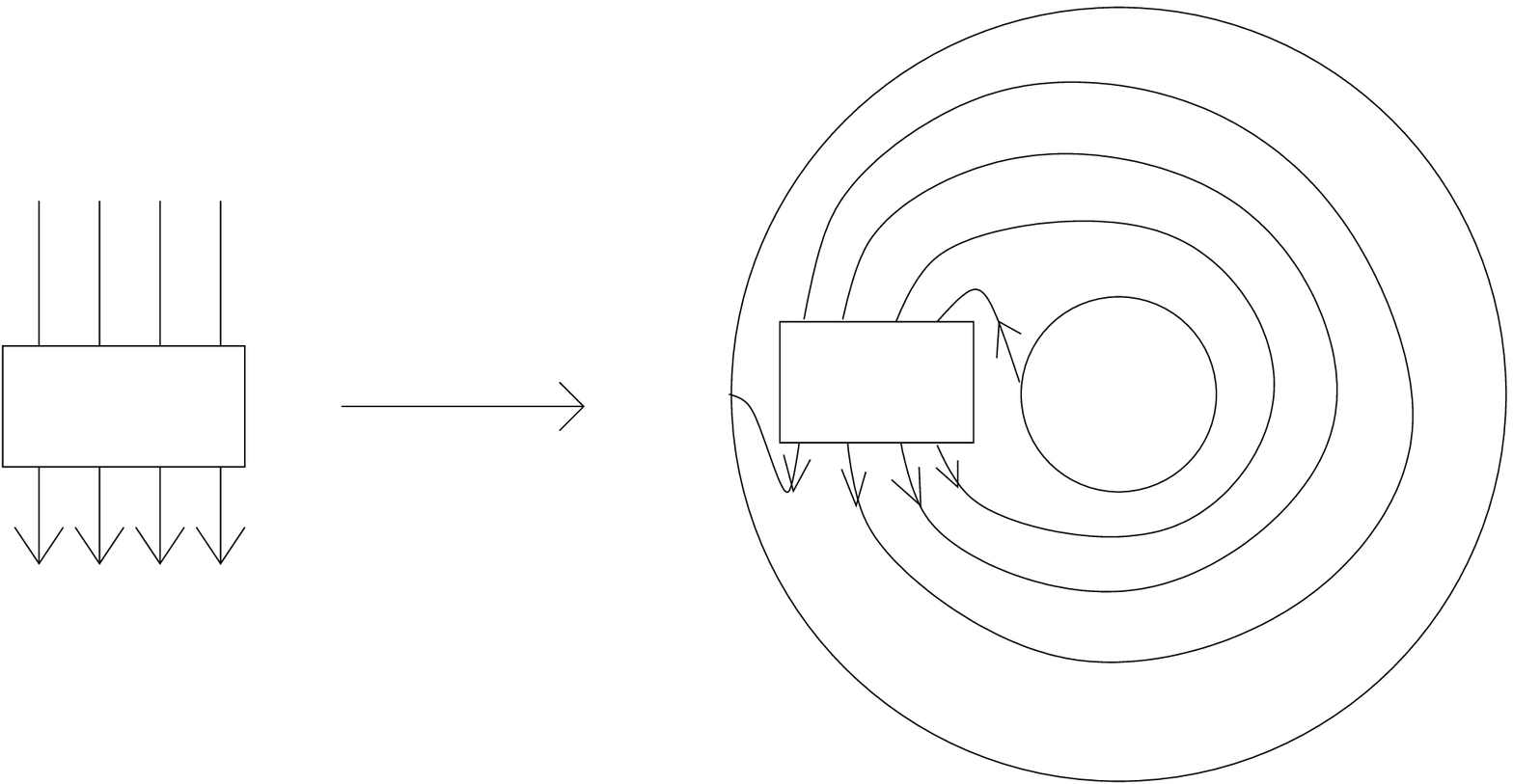,width=70\unitlength}
     \caption{\label{F-jprime}Map $\Delta'_n$ from $R_n^n$ to $\sap$.}
  \end{center}
\end{figure}

For any integer $n\geq 1$ we have a linear map $\Delta'_n:R_n^n\to \sap$
as shown in figure \ref{F-jprime}.
We abbreviate $\Delta'_n$ by $\Delta'$ if the index $n$ is obvious from the
context. 
With this notation, the diagram $a$ from figure \ref{F-arca} is the
image of the identity braid $1_2$ of $R_2^2$ under the map $\Delta'_2$.

We denote the image of $R_n^n$
under the map $\Delta'_n$ by $\sa_n$.
We define $\sap_+$ to be the submodule of $\sap$ spanned by all submodules
$\sap_0,\sap_1,\ldots$.

 The diagrams $A\cdot B$ and
$B\cdot A$ are not regularly isotopic in general but they are equal modulo
the skein relations as will be shown in lemma \ref{L-cprico}.

We note that $\sap$ can be turned into an algebra over $\sa$ in two ways.
Let $A$ be a diagram in $\sap$ and $\gamma$ be a diagram in $\sa$. Then
$\gamma\cdot A$ is defined as putting $\gamma$ above $A$, and $A\cdot \gamma$
is defined as putting $\gamma$ below $A$. Morton uses in \cite{Morton} the notation
$l(\gamma, A)$ resp. $r(A,\gamma)$ for the operation of $\sa$ on $\sap$ on the
left resp. right (up to a turn over of the annulus).
\begin{lemma}[\cite{Morton}]\label{L-cprico}
$\sap$ is commutative.
\end{lemma}
\begin{proof}
Let $D$ be a diagram in $\sap$. By induction on the 
number of crossings of $D$, we can write $D$ as a scalar linear
combination of diagrams, $D=\alpha_1 D_1+\cdots+\alpha_k D_k$,
such that the single arc of any diagram $D_i$ is totally descending and that
it lies completely below any other component of $D_i$. Hence,
$D=\gamma_1a^{i_1}+\cdots+\gamma_ka^{i_k}$ for elements $\gamma_1,
\ldots,\gamma_k$ of $\sa$ and $i_j\in\myZ$.

The commutativity of $\sap$ now follows from the observation that terms of the
kind $\gamma a^i$ for $\gamma\in\sa$ and $i\in\myZ$ commute with each other.
\end{proof}
\begin{remark}
The above proof suggests that $\sap_+$ is the polynomial ring in $a$
with coefficients in $\sa_+$. This turns out to be true. One has to check
that the powers of $a$ are linearly independent over the scalars $\sa_+$
(see \cite{Gilmer} and \cite{Lukac}).
\end{remark}

\begin{figure}
  \begin{center}
    \setlength{\unitlength}{1mm}
    \begin{picture}(0,0)
    \end{picture}
     \epsfig{file=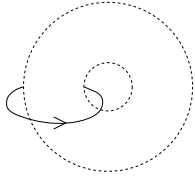,width=25\unitlength}
     \caption{\label{F-closearc}The additional arc for the closure.}
  \end{center}
\end{figure}
We define a closing operation $r\mapsto \hat{r}$
 from $\sap$ to $\sa$
which means adding the arc in figure \ref{F-closearc} from above to 
a diagram $r$.
In order that this is possible, the annulus for $\sa$ has to be slightly
larger than $\sap$.
The framing of the diagram $\hat{r}$ is defined to be its blackboard
framing. We remark that this closing operation
is not an algebra homomorphism. 

\section{Basic skein relations in $\sap$}

We recall that the quasi-idempotent $a_i\in R_i^i$ satisfies
$a_i a_i=\alpha_i a_i$ for a non-zero scalar $\alpha_i$. We define
$h_i=\frac{1}{\alpha_i}\Delta(a_i)\in \sa_i$ and $h_i'=
\frac{1}{\alpha_i}\Delta'(a_i)\in\sap_i$ for any integer
$i\geq 0$. We define $h_i=0$ and $h_i'=0$ for any integer $i<0$.

\begin{figure}
  \begin{center}
    \setlength{\unitlength}{1mm}
    \begin{picture}(0,0)
    \put(-6,16){$\frac{\qi{i+1}}{\alpha_{i+1}}$}
    \put(5,16){$a_{i+1}$}
    \put(38,16){$=\frac{1}{\alpha_i}$}
    \put(55,15){$a_i$}
    \put(85,16){$+\frac{s^{-1}\qi{i}}{\alpha_i}$}
    \put(104,16){$a_i$}
    \end{picture}
     \epsfig{file=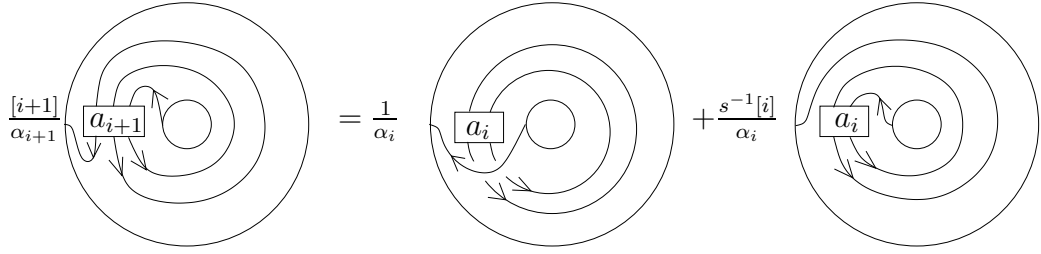,width=130\unitlength}
     \caption{\label{F-lemaiple}Depiction of lemma \ref{L-aiple}.}
  \end{center}
\end{figure}

The following lemma is a direct consequence of the decomposition
\[
  a_{i+1}
 =(a_i\tensor 1_1) (1_{i+1}+ s\sigma_i+s^2\sigma_i\sigma_{i-1}
       +\cdots+s^i\sigma_i\sigma_{i-1}\cdots\sigma_1)
\]
of the quasi-idempotent $a_{i+1}\in R_{i+1}^{i+1}$ (see \cite{Morton},
\cite{Lukac} for a proof).
The lemma is depicted in figure \ref{F-lemaiple}.

\begin{lemma}\label{L-aiple}
We have
\[
  \qi{i+1}h_{i+1}'
  =eh_i
   +s^{-1}\qi{i}h_i'a
\]
in $\sap_+$ for any integer $i\geq 0$.
\end{lemma}

We define
\[
  t_i=h_ie - eh_i\in \sap_+
\]
for any integer $i$.
We remark that $t_i=0$ for $i\leq 0$.
\begin{lemma}\label{L-yiAia}
We have
\[
  t_i=
  (s^{-i}-s^i)h_i'a
\]
for any integer $i\geq 0$.
\end{lemma}
\begin{proof}
We have
\[
   \qi{i+1}h_{i+1}'
  =eh_i+
   s^{-1}\qi{i}h_i'a
\]
by lemma \ref{L-aiple} for any integer $i$.
By applying the map $\rho$ from section \ref{sec-skeinpl} we get
\[
   \qi{i+1}h_{i+1}'
  = h_i e+s\qi{i}h_i' a
\]
because $h_i$, $h'_i$ and $h'_{i+1}$ are invariant under $\rho$
by lemma \ref{L-alpn}.
These two equations imply that $t_i=h_ie-eh_i=(s^{-i}-s^i)h'_ia$.
\end{proof}

\begin{figure}
  \begin{center}
    \setlength{\unitlength}{1mm}
    \begin{picture}(0,0)
    \put(16,22){$a_i$}
    \put(35,22){=}
    \put(41,22){$v\ \cdot$}
    \put(65,22){$a_i$}
    \put(80,22){$=\ \ \ vs^{1-i}\cdot$}
    \put(123,22){$a_i$}
    \end{picture}
     \epsfig{file=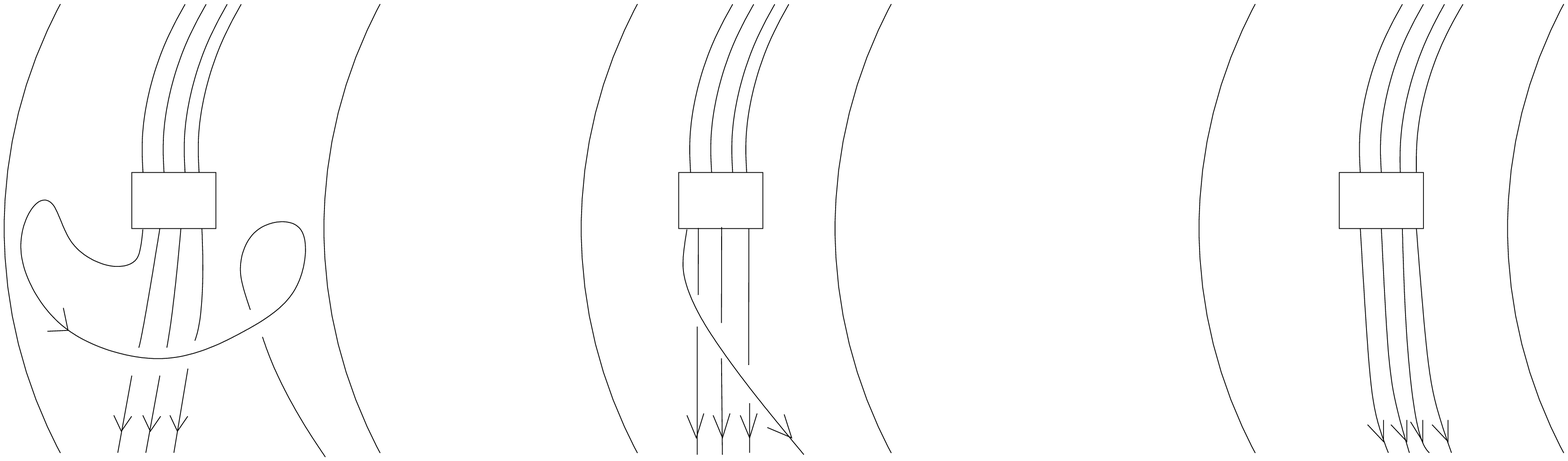,width=140\unitlength}
     \caption{\label{F-closeai-a}The closure of $\Delta'(a_i)a$.}
  \end{center}
\end{figure}

\begin{corollary}\label{C-hatyi}
We have
\[
  \hat{t}_i=(s^{1-2i}-s)vh_i
\]
for any integer $i$.
\end{corollary}
\begin{proof}
From lemma \ref{L-yiAia} and the skein relation in figure \ref{F-closeai-a}
we deduce that
\begin{eqnarray*}
  \hat{t}_i&=&(s^{-i}-s^i)vs^{1-i}h_i\\
           &=&(s^{1-2i}-s)vh_i
\end{eqnarray*}
for any integer $i\geq 0$. This equation holds for negative integers
$i$ as well because $h_i$ and $t_i$ are equal to zero for negative $i$.
\end{proof}

\begin{lemma}\label{L-yiyj}
We have
\[
  \left|
    \begin{array}{ll}
       t_i&t_{i+1}\\
       t_j&t_{j+1}
    \end{array}
  \right|
  =
  (s^2-1)
  \left|
    \begin{array}{ll}
       eh_i&t_{i+1}\\
       eh_j&t_{j+1}
    \end{array}
  \right|
\]
in $\sap$ for any integers $i$ and $j$.
\end{lemma}
\begin{proof}
If either $i$ or $j$ is negative then the lemma is obviously true
because $t_k=0$ for $k\leq 0$.
Let $i\geq 0$ and $j\geq 0$ from now on.
The combination of lemmas \ref{L-aiple} and \ref{L-yiAia} implies that
\[
  \qi{i+1}h'_{i+1}=eh_i+\frac{1}{1-s^2}t_i.
\]
We have $(s^{-1}-s)\qi{j+1}h_{j+1}'a=t_{j+1}$ by lemma \ref{L-yiAia}.
Multiplication of the above equation by these terms gives
\[
  (s^{-1}-s)\qi{i+1}\qi{j+1}h_{i+1}'h_{j+1}'a
  =
  (eh_i)t_{j+1}+\frac{1}{1-s^2}t_it_{j+1}.
\]
The left hand side of the above equation is invariant under the
interchange of $i$ and $j$ because $\sap$ is commutative,
and thus the right hand side is invariant
under this interchange. Hence,
\[
  (eh_i)t_{j+1}+\frac{1}{1-s^2}t_it_{j+1}
  =
  (eh_j)t_{i+1}+\frac{1}{1-s^2}t_jt_{i+1}
\]
which is equivalent to our claim.
\end{proof}
Lemma \ref{L-yiyj} is the stepping stone from skein calculations
to the following algebraic calculations.

\section{Determinantal calculations in $\sap$}
\begin{lemma}\label{L-ehex}
For any integer $r\geq 2$ and integers $i_1,i_2,\ldots,i_r$
we have an equality of $(r\times r)$-determinants in $\sap_+$
\[
  \left|
  \begin{array}{llll}
   h_{i_1}e&\cdots&h_{i_1+r-2}e&t_{i_1+r-1}\\
   \vdots&&\vdots&\vdots\\
   h_{i_r}e&\cdots&h_{i_r+r-2}e&t_{i_r+r-1}
  \end{array}
  \right|
  =s^{2(r-1)}
  \left|
  \begin{array}{llll}
   eh_{i_1}&\cdots&eh_{i_1+r-2}&t_{i_1+r-1}\\
   \vdots&&\vdots&\vdots\\
   eh_{i_r}&\cdots&eh_{i_r+r-2}&t_{i_r+r-1}
  \end{array}
  \right|.
\]
\end{lemma}
\begin{proof}
Since $t_i=h_ie-eh_i$, we deduce from lemma \ref{L-yiyj}
by the multilinearity of the determinant that
\begin{eqnarray}\label{E-edet2}
  \left|
    \begin{array}{ll}
       h_ie&t_{i+1}\\
       h_je&t_{j+1}
    \end{array}
  \right|
  =
  s^2
  \left|
    \begin{array}{ll}
       eh_i&t_{i+1}\\
       eh_j&t_{j+1}
    \end{array}
  \right|
\end{eqnarray}
and
\begin{eqnarray}\label{E-zweitform}
  \left|
    \begin{array}{ll}
       t_i&t_{i+1}\\
       t_j&t_{j+1}
    \end{array}
  \right|
  =  (1-s^{-2})
  \left|
    \begin{array}{ll}
       h_ie&t_{i+1}\\
       h_je&t_{j+1}
    \end{array}
  \right|.
\end{eqnarray}  
We remark that equation (\ref{E-edet2}) is our claim in the case
$r=2$.

From now on let $r\geq 3$. We see that
\[
  \left|
  \begin{array}{lllll}
    t_{i_1}&t_{i_1+1}&t_{i_1+2}&\cdots&
       t_{i_1+r-1}\\
    \vdots&\vdots&\vdots&&\vdots\\
    t_{i_r}&t_{i_r+1}&t_{i_r+2}&\cdots&
       t_{i_r+r-1}
  \end{array}
  \right|
  =
  (1-s^{-2})  
  \left|
  \begin{array}{lllll}
   h_{i_1}e&t_{i_1+1}&t_{i_1+2}&\cdots&t_{i_1+r-1}\\
   \vdots&\vdots&\vdots&&\vdots\\
   h_{i_r}e&t_{i_r+1}&t_{i_r+2}&\cdots&t_{i_r+r-1}
  \end{array}
  \right|
\]
by developing the determinant on the left hand side by the first
two columns, applying equation (\ref{E-zweitform}) to each summand,
and redeveloping the determinant.
By doing this successively for the columns 1 and 2, 2 and 3, ..., 
$(r-1)$ and $r$, we deduce that
\[
 \left|
  \begin{array}{llll}
    t_{i_1}&\cdots&t_{i_1+r-2}&t_{i_1+r-1}\\
    \vdots&&\vdots&\vdots\\
    t_{i_r}&\cdots&t_{i_r+r-2}&t_{i_r+r-1}
  \end{array}
  \right|
  =
  (1-s^{-2})^{r-1}  
  \left|
  \begin{array}{llll}
   h_{i_1}e&\cdots&h_{i_1+r-2}e&t_{i_1+r-1}\\
   \vdots&&\vdots&\vdots\\
   h_{i_r}e&\cdots&h_{i_r+r-2}e&t_{i_r+r-1}
  \end{array}
  \right|.
\] 

On the other hand, if we use the equation from lemma \ref{L-yiyj} instead
of equation (\ref{E-zweitform}) in the above argument, we get
\[
 \left|
  \begin{array}{llll}
    t_{i_1}&\cdots&t_{i_1+r-2}&t_{i_1+r-1}\\
    \vdots&&\vdots&\vdots\\
    t_{i_r}&\cdots&t_{i_r+r-2}&t_{i_r+r-1}
  \end{array}
  \right|
  =
  (s^2-1)^{r-1}  
  \left|
  \begin{array}{llll}
   eh_{i_1}&\cdots&eh_{i_1+r-2}&t_{i_1+r-1}\\
   \vdots&&\vdots&\vdots\\
   eh_{i_r}&\cdots&eh_{i_r+r-2}&t_{i_r+r-1}
  \end{array}
  \right|.
\] 
The above two equations imply that
\[
  \left|
  \begin{array}{llll}
   h_{i_1}e&\cdots&h_{i_1+r-2}e&t_{i_1+r-1}\\
   \vdots&&\vdots&\vdots\\
   h_{i_r}e&\cdots&h_{i_r+r-2}e&t_{i_r+r-1}
  \end{array}
  \right|
  =
  s^{2(r-1)}
  \left|
  \begin{array}{llll}
   eh_{i_1}&\cdots&eh_{i_1+r-2}&t_{i_1+r-1}\\
   \vdots&&\vdots&\vdots\\
   eh_{i_r}&\cdots&eh_{i_r+r-2}&t_{i_r+r-1}
  \end{array}
  \right|.
\]
\end{proof}

The ring of symmetric functions (in countably many variables)
is freely generated as a commutative
algebra by the complete symmetric functions $h_1, h_2,\ldots$.
We therefore have a well defined map from the ring of symmetric
functions to the skein of the annulus by mapping the complete symmetric
function $h_i$ to the skein element $h_i\in\sa$.
We denote the image of the Schur function $s_{\lambda}$ by $S_{\lambda}$.
The Schur function $s_{\lambda}$ can be expressed in terms of the
complete symmetric functions via the formula
$s_{\lambda}=\det(h_{\lambda_i+j-i})_{1\leq i,j\leq l(\lambda)}$
from \cite{Macdonald}.
We therefore have
\[
  S_{\lambda}=\det(h_{\lambda_i+j-i})_{1\leq i,j\leq l(\lambda)}\in \sa_n
\]
where $n=|\lambda|$.
We recall the linear map $\Gamma$ from section \ref{sec-skeinan}.
\begin{theorem}\label{T-slambdascalar}
We have $\Gamma(S_{\lambda})=q_{\lambda}S_{\lambda}$ in $\sa_+$ with
the scalar
\[
  q_{\lambda}=\frac{v^{-1}-v}{s-s^{-1}}+vs^{-1}
   \sum_{k=1}^{l(\lambda)}(s^{2(k-\lambda_k)}-s^{2k}).
\]
\end{theorem}

\begin{proof}
For any elements $\alpha$ and $\beta$
of the skein of the annulus $\sa$ we have
$(\alpha e)\cdot(\beta e)=(\alpha\beta) e$ in $\sap$
where $e$ is the identity of $\sap$.    
 Hence
\[
  S_{\lambda}e=\det(h_{\lambda_i+j-i}e)_{1\leq i,j\leq l(\lambda)}.
\]
Similarly
\[
  eS_{\lambda}=\det(eh_{\lambda_i+j-i})_{1\leq i,j\leq l(\lambda)}.
\]
We denote $l(\lambda)$ by $n$ from now on. We remark that 
the closure $(eS_{\lambda})^{\wedge}$ is equal to $S_{\lambda}$
with a disjoint circle which can be removed at the expense of
the scalar $(v^{-1}-v)/(s-s^{-1})$. The closure $(S_{\lambda}e)^{\wedge}$
is equal to $\Gamma(S_{\lambda})$.

By the multilinearity of the determinant we can write 
the difference of any two $(n\times n)$-determinants as a telescope sum
of $n$ $(n\times n)$-determinants.
\begin{eqnarray*}
 & &\hspace{-10mm} \left|
  \begin{array}{lll}
  y_{11}&\cdots&y_{1n}\\
  \vdots&&\vdots\\
  y_{n1}&\cdots&y_{nn}
  \end{array}
  \right|-
  \left|
  \begin{array}{lll}
  z_{11}&\cdots&z_{1n}\\
  \vdots&&\vdots\\
  z_{n1}&\cdots&z_{nn}
  \end{array}
  \right| =\\
 & &
  \sum_{k=1}^n \left|
                \begin{array}{lllllll}
                y_{1\,1}&\cdots&y_{1\,k-1}&(y_{1\,k}-z_{1\,k})&z_{1\,k+1}
                   &\cdots&z_{1\,n}\\
                \vdots&&\vdots&\vdots&\vdots&&\vdots\\
                y_{n\,1}&\cdots&y_{n\,k-1}&(y_{n\,k}-z_{n\,k})&z_{n\,k+1}
                   &\cdots&z_{n\,n}
                \end{array}
               \right|.
\end{eqnarray*}
Applying this formula to the determinants for $S_{\lambda}e$ and
$eS_{\lambda}$ we get
\begin{eqnarray*}
& &\hspace{-10mm} S_{\lambda}e-eS_{\lambda}= \\
& &\sum_{k=1}^n
  \left|
  \begin{array}{lllllll}
   h_{\lambda_1}e&\cdots&h_{\lambda_1+k-2}e&t_{\lambda_1+k-1}&eh_{\lambda_1+k}
    &\cdots&eh_{\lambda_1+n-1}\\
   \vdots&&\vdots&\vdots&\vdots&&\vdots\\
   h_{\lambda_n+1-n}e&\cdots&h_{\lambda_n+k-1-n}e&t_{\lambda_n+k-n}&
    eh_{\lambda_n+k+1-n}&\cdots&eh_{\lambda_n}
  \end{array}
  \right|.\\
\end{eqnarray*}
By lemma \ref{L-ehex} we deduce that  
\begin{eqnarray*}
& &\hspace{-7mm} S_{\lambda}e-eS_{\lambda}=\\
& &\hspace{-5mm} \sum_{k=1}^n s^{2(k-1)}
  \left|
  \begin{array}{lllllll}
   eh_{\lambda_1}&\cdots&eh_{\lambda_1+k-2}&t_{\lambda_1+k-1}&eh_{\lambda_1+k}
    &\cdots&eh_{\lambda_1+n-1}\\
   \vdots&&\vdots&\vdots&\vdots&&\vdots\\
   eh_{\lambda_n+1-n}&\cdots&eh_{\lambda_n+k-1-n}&t_{\lambda_n+k-n}&
    eh_{\lambda_n+k+1-n}&\cdots&eh_{\lambda_n}
  \end{array}
  \right|.
\end{eqnarray*}
The $n$ determinants appearing here are very special because each of them
is a sum of terms of the form of a $t_i$ above a product of $h_j$'s.
Therefore the closure of each determinant
is $\hat{t}_i$ above a product of $h_j$'s. Explicitly,
\begin{eqnarray*}
& &\hspace{-7mm} (S_{\lambda}e)^{\wedge}-(eS_{\lambda})^{\wedge}=\\
& &\sum_{k=1}^n s^{2(k-1)}
  \left|
  \begin{array}{lllllll}
   h_{\lambda_1}&\cdots&h_{\lambda_1+k-2}&\hat{t}_{\lambda_1+k-1}&
       h_{\lambda_1+k}&\cdots&h_{\lambda_1+n-1}\\
   \vdots&&\vdots&\vdots&\vdots&&\vdots\\
   h_{\lambda_n+1-n}&\cdots&h_{\lambda_n+k-1-n}&\hat{t}_{\lambda_n+k-n}&
       h_{\lambda_n+k+1-n}&\cdots&h_{\lambda_n}
  \end{array}
  \right|.
\end{eqnarray*}
We know by corollary
\ref{C-hatyi} that $\hat{t}_i$ is a scalar multiple of $h_i$.
Hence
\begin{eqnarray*}
& &\hspace{-10mm}(S_{\lambda}e)^{\wedge}-(eS_{\lambda})^{\wedge}=\\
& &\sum_{k=1}^n
  \left|
  \begin{array}{lllllll}
   h_{\lambda_1}&\cdots&h_{\lambda_1+k-2}&\beta_{1\, k}
       h_{\lambda_1+k-1}&
       h_{\lambda_1+k}&\cdots&h_{\lambda_1+n-1}\\
   \vdots&&\vdots&\vdots&\vdots&&\vdots\\
   h_{\lambda_n+1-n}&\cdots&h_{\lambda_n+k-1-n}&\beta_{n\, k}
       h_{\lambda_n+k-n}&
       h_{\lambda_n+k+1-n}&\cdots&h_{\lambda_n}
  \end{array}
  \right|
\end{eqnarray*}
where $\beta_{i\, k}=s^{2(k-1)}(s^{1-2(\lambda_i+k-i)}-s)v$.
We use the notation $\rho_i=s^{2i-2\lambda_i-1}v$ and $\gamma_k=-s^{2k-1}v$,
hence $\beta_{i\, k}=\rho_i+\gamma_k$.
By the multilinearity of the determinant we get
\begin{eqnarray*}
& &\hspace{-10mm}(S_{\lambda}e)^{\wedge}-(eS_{\lambda})^{\wedge}=
   (\gamma_1+\ldots+\gamma_n)S_{\lambda}+\\
& &\sum_{k=1}^n
  \left|
  \begin{array}{lllllll}
   h_{\lambda_1}&\cdots&h_{\lambda_1+k-2}&\rho_1
       h_{\lambda_1+k-1}&
       h_{\lambda_1+k}&\cdots&h_{\lambda_1+n-1}\\
   \vdots&&\vdots&\vdots&\vdots&&\vdots\\
   h_{\lambda_n+1-n}&\cdots&h_{\lambda_n+k-1-n}&\rho_n
       h_{\lambda_n+k-n}&
       h_{\lambda_n+k+1-n}&\cdots&h_{\lambda_n}
  \end{array}
  \right|.
\end{eqnarray*}
We bring the sum over the determinants in a more
appropriate form via the general formula for variables
$w_{ij}$ and $\rho_k$,
\begin{eqnarray*}
& &\hspace{-17mm}\sum_{k=1}^n
  \left|
    \begin{array}{ccccccc}
  w_{1\,1}&\cdots&w_{1\, k-1}&\rho_1 w_{1\, k}&w_{1\, k+1}&\cdots&w_{1\,n}\\
  \vdots&&\vdots&\vdots&\vdots&&\vdots\\
  w_{n\,1}&\cdots&w_{n\, k-1}&\rho_n w_{n\, k}&w_{n\, k+1}&\cdots&w_{n\,n}\\
    \end{array}
  \right|=\\
& & \hspace{50mm}
   (\rho_1+\cdots+\rho_n)
   \left|
   \begin{array}{ccc}
    w_{1\,1}&\cdots&w_{1\, n}\\
    \vdots&&\vdots\\
    w_{n\, 1}&\cdots&w_{n\, n}
    \end{array}
   \right|.
\end{eqnarray*}
Applying this formula we get
\begin{eqnarray*}
(S_{\lambda}e)^{\wedge}-(eS_{\lambda})^{\wedge}
&=&(\gamma_1+\cdots+\gamma_n)S_{\lambda}+(\rho_1+\cdots+\rho_n)
 S_{\lambda}\\
&=&(\beta_{1\, 1}+\cdots+\beta_{n\, n})S_{\lambda}.
\end{eqnarray*}
Since $(eS_{\lambda})^{\wedge}=(v^{-1}-v)/(s-s^{-1}) S_{\lambda}$, as mentioned above,
 we deduce that
$(S_{\lambda}e)^{\wedge}=q_{\lambda}S_{\lambda}$ with
\begin{eqnarray*}
  q_{\lambda}
  &=&\frac{v^{-1}-v}{s-s^{-1}}+\beta_{1\, 1}+\cdots+\beta_{n\, n}   \\
  &=&\frac{v^{-1}-v}{s-s^{-1}}
     +vs^{-1} \sum_{k=1}^n (s^{2(k-\lambda_k)}-s^{2k}).
\end{eqnarray*}
\end{proof}

\section{Proof that $S_{\lambda}=Q_{\lambda}$}
We have to introduce some notation.
Let $\lambda$ be a Young diagram.
A \emph{standard tableau} of 
$\lambda$ is a labelling of the cells of $\lambda$ with the integers
$1,2,\ldots,|\lambda|$ which is increasing in each row from left to right
and in each column from top to bottom. The Young diagram underlying a
standard tableau $t$ is denoted by $\lambda(t)$.
The number of different standard tableaux for a Young diagram $\lambda$
is denoted by $d(\lambda)$.
\begin{theorem}\label{T-slamqlam}
$S_{\lambda}$ is equal to $Q_{\lambda}$ for any Young diagram $\lambda$.
\end{theorem}
\begin{proof}
$Q_{\lambda}$ is non-zero in $\sa$.
Since the scalars $q_{\lambda}$ and $c_{\lambda}$
from theorem \ref{T-slambdascalar} and from 
lemma \ref{L-qeigen} are equal, 
we have that $S_{\lambda}$ and $Q_{\lambda}$ are eigenvectors
with the same eigenvalue under the map $\Gamma$.
Possibly, $S_{\lambda}=0$.
The set of $Q_{\lambda}$ for all Young diagrams $\lambda$
with $n$ cells is a linear basis for $\sa_n$ by lemma \ref{L-dimcn}.
Furthermore, the eigenvalues $c_{\lambda}$ are pairwise different.

Hence, we deduce that $S_{\lambda}$ is
a scalar multiple of $Q_{\lambda}$ for any Young diagram 
$\lambda$ with $n$ cells. This scalar is a rational
function in $x$, $v$ and $s$, and it is possibly equal to zero.

We denote the Young diagram consisting of a single cell by $\bigdone$.
We have that $S_{\done}=Q_{\done}$ is the single core circle of the annulus.
Hence, $S_{\done}^n=Q_{\done}^n$ are $n$ parallel copies of the core circle of the
annulus.

Macdonald states in \cite{Macdonald} example  I.3.11
the multiplication rule for the Schur function $s_{\lambda}$ with
the first power sum $p_1$ (which is equal to $s_{\done}$).
This is a special case of the Littlewood-Richardson rule for the multiplication
of Schur functions. In this case, the rule is
\[
  s_{\mu}s_{\done}=\sum_{\begin{array}{c}\scriptstyle\mu\subset\eta\\
                                         \scriptstyle|\eta|=|\mu|+1
                         \end{array}}
                      s_{\eta}.
\]
By applying this successively we get
\[
  s_{\done}^n=\sum_{|\lambda|=n}d_{\lambda}s_{\lambda}
\]
where $d_{\lambda}$ is the number of standard tableaux of $\lambda$.
Therefore,
\begin{eqnarray}\label{E-sdone}
  S_{\done}^n=\sum_{|\lambda|=n}d_{\lambda}S_{\lambda}
\end{eqnarray}
in $\sa_n$.

Blanchet describes in \cite{Blanchet} an explicit semi-simple decomposition
of $R_n^n$.
The lemmas he uses for his construction are also contained in 
\cite{AisMor} as explained in \cite{Lukac}.
Generalizing results of Wenzl \cite{Wenzl}, Blanchet constructs 
elements $\alpha_{t}\beta_{\tau}\in R_n^n$ (which are denoted $\alpha_{t\tau}$
here) for standard tableaux $t$ and $\tau$ with $\lambda(t)=\lambda(\tau)$.
These elements satisfy $\alpha_{t\tau}\alpha_{s\sigma}=
\delta_{\tau s}\alpha_{t\sigma}$ where $\delta_{\tau s}$
is the Kronecker-Delta, i.e. the elements $\alpha_{t\tau}$ multiply like
matrix-units. An important observation is that the closure
$\Delta(\alpha_{t\tau})\in \sa_+$ is equal to zero if $t\neq \tau$, and
$\Delta(\alpha_{t\tau})=Q_{\lambda(t)}$ if $t=\tau$.
Furthermore, the identity braid
$\id_n\in R_n^n$ is equal to the sum of all $\alpha_{tt}$, i.e.
\[
  \id_n=\sum_{|\lambda(t)|=n} \alpha_{tt}.
\]
We apply the closure map $\Delta:R_n^n\to \sa_n$ to the
above equation and we get that 
\begin{eqnarray}\label{E-qdone}
  Q_{\done}^n=\sum_{|\lambda|=n} d_{\lambda} Q_{\lambda}.
\end{eqnarray}

Since $S_{\done}=Q_{\done}$, 
equations (\ref{E-sdone}) and (\ref{E-qdone}) imply that
\[
  \sum_{|\lambda|=n}d_{\lambda}(Q_{\lambda}-S_{\lambda})=0.
\]
Since
any $S_{\lambda}$ differs from $Q_{\lambda}$ by a scalar,
and $\{Q_{\lambda}\sothat\lambda\mbox{ has $n$ cells}\}$
is a basis of $\sa_n$ by lemma \ref{L-dimcn},
we get that $Q_{\lambda}=S_{\lambda}$.
\end{proof}

Theorem \ref{T-slamqlam} together with the definition of $S_{\lambda}$ implies
that 
\begin{theorem}\label{T-ringhom}
The linear map induced by $s_{\lambda}\mapsto Q_{\lambda}$ from the
ring of symmetric functions to the skein of the annulus $\sa$ is a ring
homomorphism.
\end{theorem}

\begin{acknow}
The results are part of my PhD thesis. I thank my supervisor Prof. Hugh Morton
for his guidance and help. The thesis was funded by the DAAD by a
HSP-III-grant.
\end{acknow}

\stoppage

\end{document}